\documentclass{article}
\usepackage{amsmath}
\usepackage{amssymb}
\usepackage{array}
\usepackage{geometry}
\usepackage{graphicx}
\usepackage{enumerate}
\numberwithin{equation}{section}
\begin{document}
\begin{center}
{\bf\large Scientific Endeavors of A.M. Mathai:\\
An appraisal on the occasion of his Eightieth Birthday}\\
\vskip.1cm
April 2015\\
\vskip.2cm
Edited by H.J. Haubold, United Nations\\ 
\vskip.2cm
Current affiliation of A.M. Mathai:\\
\vskip.2cm
Director, Centre for Mathematical and Statistical Sciences\\
Peechi Campus, KFRI, Peechi-680653, Kerala, India\\
directorcms458@gmail.com\\
and\\
Emeritus Professor of Mathematics and Statistics\\
McGill University, Montreal, Canada, H3A 2K6\\
mathai@math.mcgill.ca\\
\end{center}

\vskip.3cm\noindent{\bf1.\hskip.3cm Early Work in Design of Experiment and Related Problems}\\
\vskip.3cm A.M. Mathai's first paper was in the area of Design of Experiment
and Analysis of Variance in Statistics. This work was done after
finishing M.A in Mathematics at University of Toronto and waiting to
register for Ph.D, during July-August 1962. This was the first
publication which appeared in the journal {\it Biometrics}
in 1965, Mathai (1965) [{\it Biometrics}, {\bf 21}(1965),
376-385]. This problem was suggested by Professor Ralph Wormleighton
of the University of Toronto. In two-way classification with
multiple observations per cell the analysis becomes complicated
due to lack of orthogonality in the design. If two factors, such as
the amount of fertilizer used and planting methods in an
agricultural experiment to study the yield of corn, are to be tried
and if the experiment is planned to replicate $n$ times, it may
happen that some observations in some replicate may get lost and as
a result, instead of $n$ observations per cell one may have $n_{ij}$
observations in the $(i,j)$th cell. When doing the analysis of the
data, for estimating the effects of fertilizers, say,
$\alpha_1,...,\alpha_p$, one has to solve a singular system of
linear equations of the type $(I-A)\alpha=G$ where $G$ is known and
$I-A$ is singular and the unknown quantity
$\alpha'=(\alpha_1,...,\alpha_p)$ is to be evaluated. Due to
singularity, one cannot write $\alpha =(I-A)^{-1}G$. This
$A=(a_{ij})$ is the incidence matrix and has the property that all
elements are positive and $\sum_{j=1}^pa_{ij}=1$ for each
$i=1,...,p$. Mathai observed that this property means that a norm of
$A$, namely $\Vert A\Vert =\max_{i}\sum_{j=1}^p|a_{ij}|=1$ and
further, since the design is taking care of a general effect, one
can impose a condition on $\alpha_1,...,\alpha_p$ such as
$\alpha_1+...+\alpha_p=0$. Now, consider $A$ being rewritten as
$A=A-C+C$ where $C$ is a matrix where all the first row elements are
equal to the median $a_1$ of the first row elements of $A$, all the
second row elements are the median $a_2$ of the second row elements
of $A$ and so on. Now, by using the conditions on $\alpha_i$'s,
$C\alpha=O$ (null). Then
$$(I-A)\alpha=G\Rightarrow (I-B)\alpha=G
$$where $B=(b_{ij}), b_{ij}=a_{ij}-a_i, j=1,...,p$ or
$\sum_{j=1}^p|b_{ij}|=\sum_{j=1}^p|a_{ij}-a_i|=$ sum of absolute deviations
from the median $a_i$, which is the least possible. Hence the norm
$\Vert B\Vert =\max_{i}\sum_{j=1}^p|a_{ij}-a_i|$ is the least possible
and evidently $<1$. Then
$$(I-B)\alpha=G\Rightarrow \alpha=(I-B)^{-1}G=(I+B+B^2+...)G$$Note
that the convergence of the matrix series is made the fastest
possible due to the fact that the mean absolute deviation is least
when taken from the median. Thus, successive approximations are
available from $BG,B^2G,...$ but for all practical purposes of
testing hypotheses it is found that the approximation $\alpha\approx
BG$ is sufficient. This approximation avoids matrix inversion or
other complicated operations except one matrix multiplication,
namely $BG$. Encouraged by this work, the thesis was written on
sampling distributions under missing values. A concept
called ``dispersion theory'' was also developed in the thesis. It is
shown that statistical decision making is nothing but a study of a
properly defined measure of scatter or dispersion in random
variables. Some dispersions are also defined in terms of some norms
or metrics. Papers were published on the
concept, in the journal {\it Metron}, {\bf
XXVII-34.1-2}(1968), 125-135.

\vskip.3cm\noindent{\bf 2.\hskip.3cm Work on Generalized Distributions} 

\vskip.3cm This work started in 1965 when  R.K.
Saxena from Jodhpur, India, was visiting McGill University as a
post-doctoral fellow of Charles Fox, the father of Fox's H-function.
Mathai's group is responsible to call this function as Fox's
H-function. Such Mellin-Barnes type representations were available
from 1888 onwards but since Charles Fox revived the whole area and
given a new life it was decided to call the function as Fox's
H-function. Mathai translated some statistical problems in terms of
special functions and Saxena immediately gave the solutions. Several general densities were
introduced. General compatible structures for conditional densities
and prior densities, so that the unconditional and posterior
densities could be easily evaluated in Bayesian analysis problems,
were investigated. The paper got
published in the Annals [Mathai and Saxena, {\it Ann. Math.
Statist.},{\bf 40}(1969), 1439-1448].\\

\vskip.2cm Mathai decided to study the
area of special functions and the result of this study is the book
from Oxford University Press; Mathai (1993): {\it A Handbook of Generalized Special
Functions for Statistical and Physical Sciences}, Oxford University Press 1993.\\

\vskip.3cm\noindent{\bf 2.1.\hskip.3cm Early work on multivariate analysis}
\vskip.3cm
Some problems from multivariate
statistical analysis were posed by  Mathai, and Saxena could give
the solutions in terms of G and H-functions. These functions were not
computable and hence difficult to utilize in statistics or mathematics. 
This prompted Mathai and Saxena to look into computable
series forms for G and H-functions and Mathai developed an operator which could solve
the difficulties and computable series forms could be obtained. The
work of Mathai and Saxena in the area of special functions
resulted in the following books : A.M. Mathai
and R.K. Saxena ({\it Generalized Hypergeometric Functions with
Applications in Statistics and Physical Sciences}, Springer-Verlag,
Heidelberg and New York, Lecture Notes No.348, 1973; {\it The
H-function with Applications in Statistics and Other Disciplines},
Wiley Eastern, New Delhi and Wiley Halsted, New York, 1978); A.M. Mathai and H.J. Haubold, {\it Special
Functions for Applied Scientists}, Springer, New York, 2008, A.M.
Mathai, R.K. Saxena and H.J. Haubold, {\it The H-function: Theory
and Applications}, Springer, New York 2010. 

\vskip.2cm When exploring
statistical distributions and their structural decompositions Mathai
established several results in characterizations of densities, see
Mathai ({\it Canadian Mathematical  Bullettin}, {\bf 9}(1966),
95-102, {\bf 10}(1967), 239-245; {\it South African Journal of
Statistics}, {\bf 1}(10)(1967), 43-48), Gordon and Mathai ({\it
Annals of Mathematical Statistics}, {\bf 43}(1972), 205-229). These works and others' related results were put together and
brought out a monograph on characterizations, see A.M. Mathai and
G. Pederzoli, {\it Characterizations of the Normal Probability
Law}, Wiley Eastern, New Delhi and Wiley Halsted, New York, 1977.

\vskip.3cm\noindent{\bf 3.\hskip.3cm Work in Multivariate Analysis}\\

\vskip.3cm Mathai had already noted the densities of several
structures could be written in terms of G and H-functions.
Consider $x_1,x_2,...,x_r,x_{r+1},...,x_k$ mutually independently
distributed positive random variables such as exponential variables,
type-1 or type-2 beta variables or gamma variables or generalized
gamma variables etc. Consider the structures
$$u=\frac{x_1...x_r}{x_{r+1}...x_k}, ~~ v=\frac{x_1^{\delta_1}...x_r^{\delta_r}}{x_{r+1}^{\delta_{r+1}}...x_k^{\delta_k}}\eqno(3.1)
$$where $\delta_1,..,\delta_k$ are some arbitrary real powers.
Then taking the Mellin transforms or the $(s-1)$th moments of $u$
and $v$ and then taking the inverse Mellin transform one can write
the density of $u$ as a G-function in most cases or as a H-function,
and that of $v$ as a H-function. Product of independently
distributed type-1 beta random variables has the same structure of
general moments of the likelihood ratio criterion or
$\lambda$-criterion, or a one-to-one function of it, in many of
testing hypotheses problems connected with one or more multivariate
Gaussian populations and exponential populations. This showed that
one could write the exact densities in the general cases as
G-functions in most of the cases. Mathai was searching for
computable representations in the general cases. 
\vskip.2cm During
one summer camp at Queen's University, Kingston, Ontario, Canada,
Mathai met P.N. Rathie, a post-doctoral fellow of L.L. Campbell of
Queen's University, again from Jodhpur, India, and also a student of
R.K. Saxena. They started the collaboration in
information theory and at the same time investigated ways and means
of putting G-function in nice computable series form. First they
developed an operator and later Mathai perfected it, see
Mathai ({\it Annals of the Institute of Statistical
Mathematics},{\bf 23}(1971), 181-197).  The operator is of the form
$$G_{\nu}=[\frac{\partial}{\partial s}+(-\ln x)]^{\nu}\eqno(3.2)
$$which is an operator operating on the integrand in the Mellin-Barnes
representation of the density functions when the densities are
written in terms of G-functions. By using this operator, general
series expansions are obtained for G-functions of the types
$G_{0,p}^{p,0}$, which is coming from product of independent
gamma variables, $G^{p,0}_{p,p}$, which is coming from product
of independent type-1 beta variables, $G_{p,p}^{p,p}$, which is
coming from product of independent type-2 beta variables and the
general $G_{p,q}^{m,n}$, see Mathai ({\it Metron}, {\bf 28}(1970),
122-146; {\it Mathematische Nachrichten}, {\bf 48}(1970), 129-139;
{\it South African Journal of Statistics}, {\bf 5}(1971), 71-90),
Mathai and Rathie ({\it Royal Belg. Akad. Class des Sci.}, {\bf 56}(1970),
1073-1084; {\it Sankhya Series A}, {\bf 33}(1983), 45-60), Mathai and Saxena
({\it Kyungpook Mathematics Journal}, {\bf 12}(1972), 61-68; Book:
{\it Generalized Hypergeometric Functions with Applications in Statistics
and Physical Sciences}, Springer Lecture Notes No. 348, Heidelberg and New York, 1973).\\

\vskip.2cm By using the same operator in (3.2) the exact distributions of
almost all $\lambda$-criteria associated with tests of hypotheses on the
parameters of one or more Gaussian populations and exponential populations
are worked out, see Mathai ({\it Publ. l'ISUP Paris}, {\bf 19}(1970), 1-15;
{\it Journal of the Indian Statistical Association}, {\bf 8}(1970), 1-17;
{\it Annals of the Institute of Statistical  Mathematics}, {\bf 23}(1971),
181-197; {\it Trabajos de Estadistica}, {\bf 23}(1972), 67-83, 111-124;
 {\it Skand. Aktuar.}, {\bf 55}(1972), 193-198; {\it Sankhya Series A},
 {\bf 34}(1972), 161-170; {\it Annals of the Institute of Statistical Mathematics},
 {\bf 24}(1972), 53-65, {\bf 25}(1973), 557-566), Mathai and Rathie
 ({\it Journal of Statistical Research}, {\bf 4}(1970), 140-159; {\it Annals of the
Institute of Statistical Mathematics}, {\bf 22}(1970), 69-116; {\it
Statistica}, {\bf 31}(1971), 673-688; {\it Sankhya Series A}, {\bf
33}(1971), 45-60; {\it Annals of Mathematical Statistics}, {\bf
42}(1971), 1010-1019). Mathai  popularized
Mellin transform techniques, and special function technique in
general, in statistical distribution theory. Exact
distributions of almost all $\lambda$-criteria, in the null and
non-null cases, are given in explicit computable forms for the most
general cases by Mathai and his co-researchers.
The exact distributions in some non-null cases could not be obtained
for the general cases. For example, in testing equality of covariance
matrices or equality of populations in $k$ multivariate normal
populations are still open problems for $k\ge 3$, in the sense that
some representations for the general case are not available.\\

\vskip.3cm\noindent{\bf 3.1.\hskip.3cm Development of 11-digit accurate percentage points for multivariate test statistics}\\
\vskip.3cm
Even after giving the explicit computable series forms for the various
exact distributions of test statistics in the null (when the hypothesis
is true) and non-null (under the alternate hypothesis) for the general
parameters, the series forms were complicated and exact percentage points could
not be computed. When Mathai visited University of Campinas in Brazil he met
the physicist R.S. Katiyar. After six months of joint
work of simplifying the complicated gamma products, psi and zeta functions,
Katiyar was able to come up with a computer program. The first
paper in the series giving the exact percentage points up to 11-digit accurate
was produced. This paper made all the complicated
theory usable in practical situations of testing of hypotheses in multivariate
 statistical analysis. The paper appeared in Biometrika and other papers followed, see Mathai and Katiyar
 ({\it Biometrika}, {\bf 66}(1979), 353-356; {\it Annals of the Institute
 of   Statistical Mathematics}, {\bf 31}(1979), 215-224; {\it Sankhya Series B},
 {\bf 42}(1980), 333-341), Mathai ({\it Journal of Statistical Computation and
 Simulation}, {\bf 9}(1979), 169-182).

\vskip.3cm\noindent{\bf 3.2.\hskip.3cm Development of a computer algorithm for nonlinear least squares}\\

\vskip.3cm After developing a computer program
for computing exact 11-digit accurate percentage points from
complicated series forms of the exact densities of
$\lambda$-criteria for almost all multivariate test statistics, the
problem of developing a computer program for non-linear least
squares was re-examined. Starting with Marquardt's methods, there
were a number of algorithms available in the literature but all
these algorithms had deficiencies. There are a few (around 11)
standard test problems to test the efficiency of a computer program.
The efficiency of a computer program is measured by checking the
following two items:  In how many test functions the computer
program fails and how many function evaluations are needed to come
up to the final solution. These are the usual two criteria used in
the field to test a new algorithm. A new algorithm for non-linear
least squares was developed by Mathai and Katiyar which did not fail
in any of the test functions and the number of function evaluations
needed was least compared to all other algorithms available in the
literature. The paper was
published in a Russian journal, see Mathai and Katiyar ({\it
Researches in Mathematical Statistics (Russian)}, {\bf
207(10)}(1993), 143-157). This paper was later
translated into English by the American Mathematical Society.\\

\vskip.3cm\noindent{\bf 3.3.\hskip.3cm Integer programming problem}\\

\vskip.3cm The usual optimization problems such as optimizing a
quadratic form or quadratic expression, subject to linear or quadratic
constraints, optimizing a linear form subject to linear (linear
programming problem) or quadratic constraints etc deal with continuous
variables. When the variables are continuous then these optimization
problems can be handled by using calculus or related techniques.
Suppose that the variables can only take integer values such as
positive integers $1,2,3,...$ then the problem becomes complicated.
Many of the standard results available when the variables are
continuous are no longer true when the variables are integer-valued.
One such problem was brought to the attention of Mathai by S. Kounias.
This was solved and a joint paper was published, see Kounias and Mathai ({\it Optimization}, {\bf 19}(1988), 123-131).\\

\vskip.3cm\noindent{\bf 4.\hskip.3cm Work on Information Theory}\\

\vskip.3cm When the exact distributions for the test statistics being worked out,
side by side the work on information theory was also progressing.
 Characterizations of information and statistical concepts were the
 ones attempted as a joint venture by Mathai and Rathie. Several
 characterization theorems were established for various information
 measures and for statistical concepts such as covariance, variance,
 correlation etc, see for example, Mathai and Rathie ({\it Sankhya Series A},
 {\bf 34}(1972), 441-442; {\it Annals of the Institute of Statistical
 Mathematics}, {\bf 24}(1972), 473-483; in the book {\it Measures of
 Information and Their Applications}, IIT Bombay, pp. 1-10, 1974; in the
 book {\it Essays in Probability and Statistics}, Shinko Tsusho, Tokyo,
 pp. 607-633, 1976. This collaboration resulted in the first book in the area
 of characterizations of information measures, A.M. Mathai and P.N.
 Rathie$^{*}$, {\it Basic Concepts in Information Theory and Statistics:
 Axiomatic Foundations and Applications}, Wiley Eastern, New Delhi and
 Wiley Halsted, New York, 1975. One of the measures discussed there is
 Havrda-Charv\'at $\alpha$-generalized entropy
$$H=\frac{\int_{-\infty}^{\infty}[f(x)]^{\alpha}{\rm d}x-1}{2^{1-\alpha}-1}\eqno(4.1)
$$where $f(x)$ is a density function. This is the continuous version.
There is also a discrete analogue. The denominator is put into the form of
 the exponent of $2$ for ready applications to binary systems. When
 $\alpha\to 1$ one has $H$ in (4.1) going to the Shannon entropy
 $S=-\int_{-\infty}^{\infty}f(x)\ln f(x){\rm d}x$ and hence (4.1) is
 called an $\alpha$-generalized entropy. There are several $\alpha$-generalized
 entropies in the literature, including the one given by Mathai.
 This (4.1) in a modified form with the denominator replaced by $1-\alpha$
 is developed later by C. Tsallis, as
 the basis for the whole area of non-extensive statistical mechanics.
 The Mathai-Rathie (1975)
 book can be considered to be the first book on characterizations. As a side
 result, as an application of functional equations, Mathai and Rathie
 solved a problem in graph theory, see
 {\it Journal of Combinatorial Theory}, {\bf 13}(1972), 83-90. Other
 applications of information theory concepts in social sciences, population
 studies etc may be seen from Kaufman and Mathai ({\it Journal of
 Multivariate Analysis}, {\bf 3}(1973), 236-242), Kaufman, Mathai, Rathie
 ({\it Sankhya Series A}(1972), 441-442), Mathai ({\it Transactions
 of the 7th Prague Conference on Information Theory}, pp. 353-357).

 \vskip.3cm\noindent {\bf 4.1.\hskip.3cm Applications to real-life problems}\\
 \vskip.3cm Applications of the concepts of information measures, `entropy'
 or the measure of `uncertainty', directed divergence (a concept of pseudo-distance),
 `affinity' or closeness between populations, concept of `distance between
 social groups' etc were applied to solve problems in social statistics,
 population studies etc. Mathai had developed a generalized measure of
 `affinity' as well as `distance between social groups'. On application side,
 dealing with applications of information theory type measures, see George and
 Mathai ({\it Canadian Studies in Population}, {\bf 2}(1975), 91-100, {\bf 7}(1980),
 1-7; {\it Journal of  Biosocial Sciences (UK)}, {\bf 6}(1975), 347-356;
 {\it The Manpower Journal}, {\bf 14}(1978), 69-78).\\

\vskip.3cm\noindent{\bf 5.\hskip.3cm Work on Biological Modeling}

\vskip.3cm During one of the visits of Mathai to the Indian
Statistical Institute in Calcutta, India, he came across the
biologist T.A. Davis. Davis had a number of problems for which he
needed answers. He had a huge collection of data on the number of
petals in certain flowers of one species of plant. He noted that the
petals were usually 4 in each flower but sometimes the number of
petals was 5. He wanted to know whether the occurrence of 5-petaled
flowers showed any pattern. His data were insufficient to come up
with any pattern. Patterns, if any, would be connected to genetical
factors. Then he had a question about how various patterns come in
nature, in the growth of leaves, flowers, arrangements of petals and
seeds in flowers etc and whether any mathematical theory could be
developed to explain these. Then he brought in the observations on
sunflower. When we look at  flowers, certain flowers such as rose
flower, sunflower etc look more beautiful than other flowers. This
appeal is due to the arrangements of petals, florets, and color
combinations. When we look at a sunflower at the florets or at the
seed formations, after the florets dry up, we see some patterns in
the arrangements of these seeds on the flower disk called capitulum.
The seeds look like arranged along some spirals coming from the
periphery going to the center. Let us call these as radial spirals.
If one marks a point on the periphery and then one looks to the left
of the mark one sees one set of radial spirals and if one looks to the
right one sees a different set of radial spirals going in the
opposite direction. The numbers of these two sets are always two
successive numbers from a Fibonacci sequence of numbers
$1,1,2,3,5,8,13,21,...$ (the sum of two successive numbers is the
next number). Another observation made is that if one looks along a
radial spiral this spiral does not go to the center but it becomes
fuzzy after a while. At that stage if one draws a concentric circle
and then look into the inside of this circle then one will see that
if one started with the pair $(13,21)$, then this has shifted to
$(8,13)$ and then to $(5,8)$ and so on. The same sort of arrangement
can be seen in pineapple, in the arrangement of leaves on a coconut
tree crown and at many other places. If one takes a coconut crown and
project onto a circle then the positions of the leaves on the crown
form a replica of the seed arrangement in a sunflower. In a coconut
crown if the oldest leaf is in a certain direction, call it 0-th
direction then the next older leaf is not the next one to the
oldest, but it is about $\theta$ degrees either to the right or to
the left and this $\theta$ is such that $\frac{\theta}{2\pi
-\theta}=\mbox{ golden ratio}=\frac{\sqrt{5}-1}{2}$. This golden
ratio also appears at many places in nature and the above $\theta
\approx 137.5^{o}$. Davis wanted a mathematical explanations for
these and related observations. These observations were made by
biologists over centuries. Many theories were also available on the
subject. All the
theories were trying to explain the appearance of radial spirals.
Mathematicians try with differential equations and others from other
fields try with their own tools. Mathai figured out that the radial
spirals that one sees may be aftermath of something else and radial
spirals are not generated per se. Also the philosophy is that nature
must be working on very simple principles. If one buys sunflower seeds from a shop or look at
sunflower seeds on a capitulum the seeds are all of the same
dimensions if one takes one from the periphery or from any other spot
on the capitulum. Such a growth can happen if something is growing
along an Archimedes' spiral, which has the equation in polar
coordinates $r=k\theta$ after one leaves the center. Davis' artist
was asked to mark points on an Archimedes' spiral, differentiating
from point to point at $\theta=\approx 137.5^{o}$, something like a
point moving along Archimedes' spiral at a constant speed so that
when the first points reaches $\theta$ mark a second point starts,
both move at the same speed whatever be the speed. When the second
point comes to the mark $\theta$ a third point starts, and so on.
After creating a certain number of points, may be 200 points, remove
the Archimedes' spiral from the paper and fill up the space with any
symmetrical object, such as circle, diamonds etc,  with those points
being the centers. Then if one looks from the periphery the two types
of radial spirals can be seen. No such spirals are there but it is
one's vision that is creating the radial spirals. Thus a sunflower
pattern was recreated from this theory and Mathai and Davis proposed
a theory of growth and forms. Consider a capillary a very thin tube
with built-in chambers. Consider a viscous fluid being continuously
pumped in from the bottom. The liquid enters the first chamber. When
a certain pressure is built up, an in-built valve opens and the
fluid moves into the second chamber and so on. Suppose that
the tube opens in the center part of a pan (with a hole at the center).
If the pan is fully sealed so that the only force acting on the
liquid is Earth's gravity. The flow of the liquid will be governed
by the functional equation $f(\theta_1)+f(\theta_2)=f(\theta_1+\theta_2)$
whose continuous solution is the linear function $f(\theta)=k\theta$.
This is Archimedes' spiral. \\
\vskip.2cm The paper was sent for publication in the journal of
Mathematical Biosciences the editor `enthusiastically
accepted for publication'. In this paper, Mathai and Davis$^{*}$({\it Mathematical
Biosciences}, {\bf 20}(1974), 117-133), a theory of growth and form is proposed.
This theory still stands and since then there were many papers in physics,
chemistry and other areas supporting various aspects of the theory and
none has disputed the theory so far. In 1976 the journal has taken
Mathai-Davis sunflower head as the cover design for the journal and it is still the cover design.\\

\vskip.3cm\noindent{\bf 5.1.\hskip.3cm Work on coconut tree crown}\\

\vskip.3cm The coconut crown was also examined from many mathematical points
of view and found to be an ideal crown. This paper may be seen from Mathai
and Davis ({\it Proceedings of the National Academy of Sciences, India},
{\bf 39}(1973), 160-169).\\

\vskip.3cm\noindent{\bf 5.2.\hskip.3cm Engineering wonder of Bayya bird's
nest and other biological problems}\\

\vskip.3cm Further problems looked into by Mathai and Davis are the
following: (1) The engineering aspect of the egg chamber of bayya
bird's nest. The nest hangs from the tips of tree branches, the
mother bird goes into the egg chamber through the tail opening of
the nest, the nest oscillates violently during heavy winds or storms
but no egg comes out of the egg chamber and fall through the tail
opening. Naturally the tail opening is bigger than the diameter of
the eggs because the mother bird goes through that opening. This
shape, beng an engineering marvel, was examined by Mathai and Davis.
(2) thermometer
birds in Andaman Nicobar Islands; (3) transfer of Canadian Maple Syrup
technology in the production of palm sugar and jaggery in Tamilnadu,
India; (4) Nipa palms to prevent sea erosion along Kannyakumari sea coast;
(5) rejuvenation of Western Ghats in Kannyakumari region.
All these projects were undertaken jointly by the Centre for
Mathematical Sciences, Trivandrum Campus (CMS) where A.M. Mathai was
the Honorary Director and Haldane Research Institute of Nagarcoil, Tamilnadu (HRI)
where T.A. Davis was the Director and A.M. Mathai was the Honorary Chairman.
Earlier to these studies, George
and Mathai had done work in population problems, especially in the
study of inter-live-birth intervals, that is, the interval between two
live births among women in child-bearing age group, see George and
Mathai ({\it Sankhya Series B}, {\bf 37}(1975), 332-342; {\it Demography
of India}, {\bf 5}(1976), 163-180; {\it The Manpower Journal}, {\bf 14}(1978),
69-78). Here, Mathai had introduced the concepts of affinity and distance between social groups.\\

\vskip.3cm \noindent{\bf 5.3.\hskip.3cm Introducing the phrase
`statistical sciences'}\\
\vskip.3cm By 1970 Mathai was working to establish a Canadian
statistical society and a Canadian journal of statistics. The phrase
`statistical sciences' was framed and defined it as a systematic and
scientific study of random phenomena so that the theoretical
developments of probability and statistics and applications in all
branches of knowledge will come under the heading `statistical
sciences', and random variables as an extension of mathematical
variables or mathematical variables as degenerate random variables.
After launching Statistical Science
Association of Canada, the term `statistical science' became a
standard phrase. Journals and organizations started using the name
`statistical science'. Mathai was responsible to introduce these
terms into scientific literature. \vskip.2cm When G.P.H. Styan, a
colleague of Mathai, was editing the news bulletin of the Institute
of Mathematical Statistics he posed the question whether the phrase
`statistical science' was ever used before launching statistical
science association of Canada. There was a response from a Japanese
scientist claiming that he had used the term `statistical science'
before. Incidentally, later the Institute of Mathematical Statistics
changed the name of Annals of Mathematical Statistics to Annals of
Statistics and hence that name was no longer available when
statistical science association of Canada changed its name back to
the original proposed name Statistical Society of Canada.

\vskip.3cm\noindent{\bf 6.\hskip.3cm Work on Probability and Geometrical Probabilities}\\

\vskip.3cm Work in mathematical statistics and special functions continued.
As a continuation of the investigation of structural properties of densities,
Mathai came across the distributions of lengths, areas and volume contents of
random geometrical configurations such as random distance, random area,
random volume and random hyper-volume. All the theories of G and H-functions,
products and ratios of positive random variables etc could be used in
examining the distributional aspects of volume of random parallelotopes
and simplices. By analyzing the structure of general moments, Mathai noted
that these could be generated by products of independently (1) gamma
distributed points, (2) uniformly distributed points, (3) type-1 beta
distributed points, (4) type-2 beta distributed points. Out of these, (1)
fell into the category of $G_{0,p}^{p,0}$, the second and third fell into
$G_{p,p}^{p,0}$ category and (4) fell into $G_{p,p}^{p,p}$ category, for all
of which the necessary theory was already developed by Mathai and his team.
Papers were published on the distributional aspects, see Mathai ({\it Sankhya
Series A}. {\bf 45}(1983), 313-323; Mathai and Tracy ({\it Communications in
Statistics A}, {\bf 12(15)}(1983), 625-632, 1727-1736; Mathai ({\it Proceedings
of ISPS VI Annual Conference}, pp. 3-8, 1987; {\it International Journal of
Mathematical and Statistical Sciences}, {\bf 3(1)}(1994), 79-109, {\bf 7(1)}
(1998), 77-96; Rendiconti del Circolo Matematico di Palermo, Serie II, Suppl.,
{\bf 65}(2000), 219-232),  Mathai and Pederzoli ({\it American Journal of
Mathematical and Management Sciences} {\bf 9}(1989), 113-139; {\it Rendiconti
del Circolo Matematico di Palermo, Serie II, Suppl.}, {\bf 50}(1997), 235-258).\\

\vskip.3cm\noindent{\bf 6.1.\hskip.3cm A conjecture in geometric probabilities}\\

 \vskip.3cm Then Mathai came across a conjecture posed by an Australian scientist
 R.E. Miles, regarding the asymptotic normality of a certain random volume coming
 from uniformly distributed random points. This was proved to be true by H.
 Ruben. In fact Ruben
 brought this area to the attention of Mathai. The structure of the random geometric configuration was known to
 Mathai and that it was a G-function of the type $G_{p,p}^{p,0}$ and Mathai
 realized that a very simple proof of the conjecture could be given by using
 the asymptotic formula, or Stirling's formula which is the first approximation
 there, for gamma functions. This was worked out and shown that the conjecture
 could be proved very easily. This paper appeared in the journal in probability,
 see Mathai ({\it Annals of Probability}. {\bf 10}(1982), 247-251). Incidentally,
 there is a mistake there. Final representation is given in terms of a confluent
 hypergeometric function ${_1F_1}$ there but it should be a Gauss hypergeometric
 function ${_2F_1}$, one parameter is missed there in writing the final form. Then
 Mathai noted that the same conjecture can be formulated in terms of type-1 beta
 distributed random points and similar conjectures could be formulated for type-2
 beta distributed random points and gamma distributed random points. These conjectures
 were formulated and solved, see Mathai ({\it Sankhya Series A}, {\bf 45}(1983), 313-323;
 {\it American Journal of Mathematical and Management Sciences}, {\bf 9}(1989),
 113-139); Mathai and Tracy ({\it Communications in Statistics A}, {\bf 12(15)}(1983),
 1727-1736; {\it Metron}, {\bf 44}(1986), 101-110).\\
 
\vskip.3cm\noindent{\bf 6.2.\hskip.3cm Random volumes and Jacobians of matrix transformations}\\

\vskip.3cm Side by side Mathai was developing functions of matrix argument.
The work in this area will be given later but its connection to geometrical
probabilities will be mentioned here. The area of stochastic geometry or
geometrical probabilities is a fusion of geometry and measure theory. When
measure theory is mixed with geometry the standard axiomatic definition for
probability measure is not sufficient. It is quite evident to see that an
additional property of invariance is needed because a geometrical object
can be moved around in a plane or in space and the probability statements
must remain the same. The famous Betrand's paradoxes or Russell's paradoxes
come from lack of invariance conditions there. The details are discussed in
the book, A.M. Mathai,{\it Introduction to
Geometrical Probability: Distributional Aspects and Applications}, Gordon
and Breach, New York, 1999. Consider a circle of radius $r$. Take two points
$A$ and $B$ at random and independently on the circumference of this circle.
Here, `at random' could mean that the probability of finding a point, such
as $A$, in an interval of length $\delta$ is $\frac{\delta}{2\pi r}$. Consider
the chord $AB$. Then $AB$ is a random chord. Let $P$ be the mid point of this
chord and $O$ the center of the circle. Then $OP$ is fixed when $AB$ is fixed
and $OP$ is perpendicular to $AB$. Consider another situation of selecting a
point $P$ at random inside the circle. This can be done by assigning probability
of finding $P$ in a region $R$ inside the circle is $\frac{R}{\pi r^2}$. If $P$
is fixed and if $P$ is the midpoint of a chord then the chord is automatically
fixed. In many ways one can geometrically uniquely determine a chord. The chord
can be made `random' by assigning probabilities in many ways. Two ways are described
above. If one asks a question, what is the probability that the length of this random
chord is less than a specified number? The answer will be different for different
ways of assigning probabilities. This is the paradox. Note that all steps in the
derivations of the answers will be correct and valid steps as per the usual axioms of probability.

\vskip.2cm In stochastic probability area the methods used are the methods from
differential and integral geometry and usually very difficult. Even if one wishes to
talk about the distribution of random volume of a parallelotope through differential
or integral geometry the process is very involved. Mathai noted that such problems
could be easily answered through Jacobians of matrix transformations. A paper was
published in advances in
applied probability, see Mathai {\it Advances in Applied Probability}, {\bf 31(2)}(1999),
 487-506). More
 papers were published, see Mathai
 ({\it Rendiconti del Circolo Matematico di Palermo, Serie II, Suppl.}, {\bf 65}(2000),
 219-232; in the book {\it Probability and Statistical Models with Applications}, pp. 293-316,
 Chapman and Hall, 2001, {\it Rendiconti del Circolo Matematico di Palermo} {\bf XLVIII}(1999),
  487-506); Mathai and Moschopoulos ({\it Statistica}, {\bf LIX(1)}(1999), 61-81; {\it
  Rendiconti del Circolo Matematico di Palermo}, {\bf XLVIII}(1999), 163-190).

\vskip.3cm\noindent{\bf 6.3.\hskip.3cm Applications in transportation problems}\\

\vskip.3cm As an application of geometrical probability problems
Mathai explored the travel distance from the suburb to city core for
circular and rectangular grid cities. Many of the European cities
are designed with a city center and circular and radial streets from
the center whereas in North America most of the cities are designed
in rectangular grids. Travel distances, time taken and associated
expenses are random quantities and related to the nature of city
design. Some problems of this type were analyzed by Mathai ({\it
Environmetrics}, {\bf 9}(1998), 617-628); Mathai and Moschopoulos
({\it Environmetrics}, {\bf  10}(1999), 791-802).

\vskip.3cm\noindent{\bf 7.\hskip.3cm Work in Astrophysics}\\

\vskip.3cm After publishing the two books on generalized
hypergeometric functions in 1973 and H-function in 1978, physicists were interested to use those results in their works.
A number of people from different parts of Germany were using these
results. The German group working in astrophysics problems were
trying to solve some problems connected with reaction rate theory.
Then H.J. Haubold, came to McGill University  with open problems where
help from special function theory was needed. After converting their
problems in terms of integral equations, Mathai noted that the basic
integral to be evaluated was of the following form:
$$I(\gamma, a,b)=\int_0^{\infty}x^{\gamma}{\rm e}^{-ax-bx^{-\frac{1}{2}}}{\rm d}x\eqno(7.1)
$$and generalizations of this integral. Note that if $a$ or $b$ is
zero then the integral can be evaluated by using a gamma integral.
Mathematically, if the nonlinear  exponent is of the form
$x^{-\frac{1}{2}}$ or of the form $x^{-\rho}, \rho>0$ it would not
make any difference. Mathai could not find any such integrals in any
of the books of tables of integrals. He noted that the integrand
consisted of integrable functions and therefore one could make
statistical densities out of them. For example,
$f_1(x)=c_1~x^{\gamma}{\rm e}^{-ax},0\le x<\infty$ is a density
where $c_1$ is the normalizing constant. Similarly $f_2(x)=c_2{\rm
e}^{-x^{\rho}},\rho>0,0\le x<\infty$ is a density where $c_2$ is the
normalizing constant. Then the structure in (7.1) can be written as
follows:
$$g(u)=\int_v\frac{1}{v}f_1(v)f_2(\frac{u}{v}){\rm d}v\eqno(7.2)
$$where $g(u)$ can represent the density of $u=x_1x_2$ where $x_1$ and $x_2$ are
independently distributed positive real scalar random variables with the densities
$f_1(x_1)$ and $f_2(x_2)$ respectively. Once the structure in (7.1) is identified as
that in (7.2) then, since the density being unique, it is only a matter of finding
the density $g(u)$ by using some other means. We can easily use the properties of
arbitrary moments. For example
$$E(u)^{s-1}=E(x_1x_2)^{s-1}=E(x_1^{s-1})E(x_2^{s-1})
$$due to statistical independence of $x_1$ and $x_2$, where $E$ denotes the expected
value. Note that $E(x_1^{s-1})$ is available from $f_1(x_1)$ and $E(x_2^{s-1})$ from
$f_2(x_2)$. Then $g(u)$ is available from the inverse, that is,
$$g(u)=\frac{1}{2\pi i}\int_{c-i\infty}^{c+i\infty}E(u^{s-1})u^{-s}{\rm d}s\eqno(7.3)
$$where $i=\sqrt{-1}$ and $c$ is determined from the poles of $E(u^{s-1})$. Thus, by
using statistical techniques the integral in (7.1) was evaluated.
After working out many results it was realized that one could also
use Mellin convolution of a product to solve integrals of the type
in (7.1). This was not seen when the method through statistical
distribution theory was devised. Various types of thermonuclear reactions,
resonant, non-resonant, depleted case, high energy
cut off case etc were investigated. The work also went into
exploring exact analytic solar models,
gravitational instability problems, solar neutrino problems,
reaction-rates, nuclear energy generation etc. The work until 1988 was summarized in the monograph Mathai and
Haubold ({\it Modern Problems in Nuclear and Neutrino
Astrophysics}, Akademie-Verlag, Berlin, 1988). Since then a lot of
work was done, some of them are the following: Haubold and Mathai
({\it Annalen der Physik}, {\bf 44}(1987), 103-116; {\it
Astronomische Nachrichten}, {\bf 308(5)}(1987), 313-318; {\it
Journal of Mathematical Physics}, {\bf 29(9)}(1988), 2069-2077; {\it
Astronomy and Astrophysics}, {\bf 203}(1988), 211-216; {\it
Astronomische Nachrichten}, {\bf 312(1)}(1991), 1-6; {\it
Astrophysics and Space Science}, {\bf 176}(1991), 51-60, {\bf
197}(1992), 153-161,{\bf 214}, 49-70,139-149, {\bf 228}(1995),
77-86, {\bf 258}(1988), 185-199; {\it American Institute of Physics,
Conference Proceedings}, {\bf 320}(1994), 102-116, {\bf 320}(1994),
89-101; {\it SIAM Review}, {\bf 40(4)}(1998), 995-997). The
collaboration also resulted in two encyclopedia articles, see
Haubold and Mathai (Sun, {\it Enclyclopaedia of Planetary
Sciences}, pp. 786-794, 1997, Structure of the Universe, {\it
Encyclopedia of Applied Physics}, {\bf 23}(1998), pp. 47-51).\\

\vskip.3cm\noindent{\bf 7.1.\hskip.3cm New results in mathematics through statistical techniques}\\
 \vskip.3cm After evaluating the basic integrals in physics problems
 by using statistical techniques, it was realized that such statistical
 techniques could be used to obtain results in mathematics. Some
 summation formulae, computable series representations, extensions
 of several mathematical identities etc were obtained through
 statistical techniques, see Mathai and Tracy ({\it Metron},
 {\bf XLII-N1-2}(1985), 117-126), Mathai and Pederzoli ({\it Metron},
 {\bf XLIII-N3-4}(1985), 157-166, Mathai and Provost ({\it Statistical Methods},
 {\bf 4(2)}(2002), 75-98).\\

 \vskip.3cm\noindent{\bf 8.\hskip.3cm Work on Differential Equations}\\

 \vskip.3cm One of the problems investigated in connection with problems
 in astrophysics was the gravitational instability problem. The problem
 was brought to the attention of Mathai by Haubold. Papers by Russian
 researchers were there on the problem of mixing two types of cosmic dusts.
 Mathai looked at it and found that by making a transformation in the
 dependent variable and by changing the operator to $t\frac{{\rm d}}{{\rm d}t}$
 instead of the integer order differential operator $D=\frac{{\rm d}}{{\rm d}t}$
 one could identify the differential equation as a particular case of the differential
 equation satisfied by a G-function. Then G-function theory could be used to
 solve the problem of mixing $k$ different cosmic dusts. Thus the first paper
 in integer order differential equation was written and published in the MIT
 journal, see Mathai ({\it Studies in Applied Mathematics}, {\bf 80}(1989), 75-93).
 Two follow-up papers were written developing the differential equation and applying
 to physics problems, see Haubold and Mathai ({\it Astronomische Nachrichten},
 {\bf 312(1)}(1991), 1-6; {\it Astrophysics and Space Science}, {\bf 214(1\&2)}(1994), 139-149).

\vskip.3cm\noindent{\bf 9.\hskip.3cm The Idea of Laplacianness of Bilinear Forms and Work on Quadratic and Bilinear Forms}\\

\vskip.3cm In the 1980's two students of Mathai, S.B. Provost
and D. Morin-Wahhab, finished their Ph.Ds in the area of
quadratic form. Mathai has also published a number of papers on
quadratic and bilinear forms by this time. Then it was decided to
bring out a book on quadratic forms in random variables. On the
mathematical side, there were books on quadratic forms but there was
none in the area of quadratic forms in random variables. Only real random variables and samples coming from Gaussian
population were considered. Later in 2005 Mathai extended the
theory to cover very general classes of populations. This aspect
will be considered later when pathway models are discussed. Only when I. Olkin pointed out 
to Mathai about the many applications of complex Gaussian case in
communication theory, after the book appeared in print, Mathai and
Provost realized that an equal amount of material was missed out: A.M. Mathai and S.B. Provost, {\it Quadratic
Forms in Random Variables: Theory and Applications}, Marcel Dekker,
New York, 1992. Work on quadratic forms and related topics may be seen
from Mathai({\it Communications in Statistics A}, {\bf 20(10)}(1991)
3159-3174; {\it International Journal of Mathematical and
Statistical Sciences}, {\bf 1(1)}(1992), 5-20; {\it Journal of
Multivariate Analysis}, {\bf 41(2)}(1992), 178-193; {\it Annals of
the Institute of Statistical Mathematics}, {\bf 44}(1992), 769-779;
{\it Journal of Applied Statistical Sciences}, {\bf 1(2)}(1993),
169-178; {\it The Canadian Journal of Statistics}, {\bf
21(3)}(1993), 277-283; {\it Journal of Multivariate Analysis}, {\bf
45}(1993), 239-246; {\it Journal of Statistical Research}, {\bf
27(1\&2)}(1993), 57-80).\\

\vskip.3cm\noindent{\bf 9.1.\hskip.3cm Chisquaredness of quadratic
forms and Laplacianness of bilinear forms}\\

\vskip.3cm Consider the following quadratic form and bilinear form:

$$Q=X'AX, X'=(x_1,...,x_p), A=(a_{ij})=A', p\times p,\hskip.5cm S=X'BY, Y'=(y_1,...,y_q), B,p\times q
$$where $x_1,...,x_p,y_1,...,y_q$ are real scalar random variables,
$A$ is a $p\times p$ matrix and $B$ is a $p\times q$ matrix, where
$p\le q$ or $p\ge q$. When $X$ is distributed as a $N_p(O,I)$ or a
$p$-variate Gaussian or normal population with mean value null and
covariance matrix an identity matrix then there is a theorem which
says that $Q$ is distributed as a chisquare with $r$ degrees of
freedom if and only if $A$ is idempotent and of rank $r$. This
result is frequently used, especially in design of
experiments and analysis of variance problems. In fact, this result
and its companion result on the independence of two quadratic forms
are the backbones of the areas of analysis of variance, analysis of
covariance, regression, model building and many others.  What is a
concept corresponding to chisquaredness of quadratic form in the
bilinear form case? It was shown by Mathai that the concept is
Laplacianness or the corresponding distribution is Laplace density
instead of chisquare density, see Mathai ({\it Journal of
Multivariate Analysis}, {\bf 45}(1993), 239-246). Apart from
introducing the concept of Laplacianness, this paper also throws light on covariance structures. When
Mathai was taking his M.A. degree in mathematics, one of the
professors in a course on multivariate analysis asked a
simple-looking question in 1962. If one has a simple random sample
from a bivariate real normal population $N_2(\mu,\Sigma), \Sigma>0$
(positive definite; standard notation), consider the sample
correlation coefficient, denoted by $r$, where

$$r=\frac{\sum_{j=1}^n(x_{ij}-\bar{x})(y_{ij}-\bar{y})}{[\sum_{j=1}^n(x_{ij}-\bar{x})^2
\sum_{j=1}^n(y_{ij}-\bar{y})^2]^{\frac{1}{2}}}.\eqno(9.1)
$$The question was what is the density of the sample covariance
$\sum_{j=1}^n(x_{ij}-\bar{x})(y_{ij}-\bar{y})/n$? The density of $r$
in the bivariate normal case, and the corresponding density for the
sample multiple correlation in the multivariate case, were already
available in the literature. The answer looked trivial because the
sample covariance is directly connected to the sample correlation.
Nobody had the answer including the professor who posed the question.
In 1990-19991 when Mathai was writing on Laplacianness, he realized
that covariance structure is nothing but a bilinear form and hence
the density of the sample covariance must be available from that of
the bilinear form. Thus, the 1962 question was answered in the
above-mentioned 1993 paper. The corresponding matrix-variate case
should also be available but nobody has worked out yet.

\vskip.3cm\noindent{\bf 9.2.\hskip.3cm Bilinear form book}\\

\vskip.3cm After publishing the quadratic form book in 1992, a lot
of work had been done on bilinear forms. Even though a bilinear form
can be written as a quadratic form, there are many properties
enjoyed by bilinear form and not enjoyed by quadratic forms.
Quadratic forms do not have covariance structures. Then T. 
Hayakawa of Japan contacted Mathai asking why not bring out a book
on bilinear form, parallel to the one on quadratic form including chapters on zonal polynomials. This book on
bilinear forms and zonal polynomials was brought out in 1995: A.M.
Mathai, S.B. Provost and T. Hayakawa, {\it Bilinear Forms
and Zonal Polynomials}, Springer, New York, 1995, in the lecture
notes series. Additional papers may be seen from Mathai and Pederzoli
({\it Journal of the Indian Statistical Society}, {\bf 3}(1995),
345-356; {\it Statistica}, {\bf LVI(4)}(1996), 4-7-41).

\vskip.3cm\noindent{\bf 10.\hskip.3cm Functions of Matrix Argument}\\

\vskip.3cm Meanwhile Mathai's work on functions of matrix argument
was progressing. These are real-valued scalar functions where the
argument is a real or complex matrix. The theory is well developed
when the argument matrix is real positive definite or hermitian
positive definite. Note that when $A$ is a square or rectangular
matrix we do not have a concept corresponding to the square root of
a scalar quantity uniquely defined. But if the matrix $A$ is real
positive definite or hermitian positive definite, written as $A>O$,
operations such as square root can be uniquely defined. Hence the
theory is developed basically for real positive definite or
hermitian positive definite matrices. Gordon and Mathai tried to
develop a matrix series and a pseudo analytic function involving
general matrices, the attempt was not fully successful but some
characterization theorems for multivariate normal population could
be established, see Gordon and Mathai ({\it Annals of Mathematical
Statistics}, {\bf 43}(1972), 205-229). Gordon has two more papers in the area, one in the {\it Annals of Statistics} and the
other in the {\it Annals of the Institute of Statistical
mathematics}. Hence the theory of real-valued scalar functions of
matrix argument is developed when the matrix is real or hermitian
positive definite. There are three approaches available in the
literature. One is through matrix-variate Laplace transform and
inverse Laplace transform developed by C. Herz and others, see for
example, Herz ({\it Annals of Mathematics},{\bf 61(3)}(1955),
474-523). Here one basic assumption is functional commutativity
$f(AB)=f(BA)$ even if $AB\ne BA$, where $A$ and $B$ are $p\times p$
matrices. Under functional commutativity we have the following
result, observing that when $A$ is symmetric there exists and
orthonormal matrix $P,PP'=I,P'P=I$ such that $P'AP=D$ where $D$ is a
diagonal matrix with the diagonal elements being the eigenvalue of
$A$. Then

$$f(A)=f(AI)=f(APP')=F(P'AP)=f(D).
$$Thus, the original function of $p(p+1)/2$ real scalar variables,
can be reduced to a function of $p$ variables, the eigenvalues of $A$.
Another approach is through zonal polynomials, developed by Constantine,
 James and others, see for example James ({\it Annals of Mathematics},
 {\bf 74}(1961), 456-469) and Constantine ({\it Annals of Mathematical
 Statistics}, {\bf34}(1963), 1270-1285). In this definition a general
 hypergeometric function with $r$ upper parameters and $s$ lower
 parameters is defined as follows:
\begin{align*}
{_rF_s}(X)&={_rF_s}(a_1,...,a_r;b_1,...,b_s; X)\\
&=\sum_{k=0}^{\infty}\sum_K\frac{(a_1)_K...(a_r)_K}{(b_1)_K...(b_s)_K}\frac{C_K(X)}{k!}&(10.1)\end{align*}where
$C_K(X)$ is zonal polynomial of order $k$,
$K=(k_1,...,k_p),k_1+...+k_p=k$, and for example,
$$(a)_K=\prod_{j=1}^p(a-\frac{(j-1)}{2})_{k_j}, (b)_k=b(b+1)...(b+k-1), (b)_0=1, b\ne 0\eqno(10.2)
$$Here also functional commutativity is assumed. They claim uniqueness
for the above series by claiming that (10.1) satisfies both the integral
equations defining matrix-variate function through the definition of
Laplace and inverse Laplace pair. The third approach is due to Mathai
and it is defined in terms of a general matrix transform or
M-transform. The M-transform of $f(-X)$ defined by the equation

$$g(\rho)=\int_{X>O}|X|^{\rho-\frac{p+1}{2}}f(-X){\rm d}X,\Re(\rho)>\frac{p-1}{2}\eqno(10.3)
$$where $\Re(\cdot)$ means the real part of $(\cdot)$. Under
functional commutativity, $f(-X)$ in (10.3) reduces to a function of
$p$ variables, the eigenvalues of $X$. But, still the left side of
(10.3) is a function of only one variable $\rho$. Hence unique
determination of $f$ through $g(\rho)$ need not be expected. It is
conjectured that $f$ is unique when $f$ is analytic in the cone of
positive definite matrices. Right now, $f(-X)$ in (10.3) remains as
a class of functions satisfying the integral equation (10.3). In
this definition, a general hypergeometric function with $r$ upper
and $s$ lower parameters will be defined as that class of functions
for which the M-transform is the following:

$$g(\rho)=\frac{\Gamma_p(a_1-\rho)...\Gamma_p(a_r-\rho)}{\Gamma_p(b_1-\rho)...\Gamma_p(b_s-\rho)},\Re(\rho)>\frac{p-1}{2}\eqno(10.4)
$$where $\Gamma_p(a)$ is the real matrix-variate gamma given by
$$\Gamma_p(a)=\pi^{\frac{p(p-1)}{4}}\Gamma(a)\Gamma(a-\frac{1}{2})...\Gamma(a-\frac{(p-1)}{2}),\Re(a)>\frac{p-1}{2}.\eqno(10.5)
$$Then that class of function $f(-X)$ is given by the equation (10.5).
It is seen that M-transform technique is the most powerful in
extending univariate results to matrix-variate cases. Some of the
results may be seen from Mathai ({\it Mathematische Nachrichten}
{\bf 84}(1978), 171-177; {\it Communications in Statistics A}, {\bf
A8(1)}(1979), 47-55, {\bf A9(8)}(1980), 795-801;{\it Annals of the
Institute of Statistical Mathematics}, {\bf 33}(1981), 35-43, {\bf
34}(1982), 591-597; {\it Sankhya Series A}, {\bf 45}(1983),
313-323; {\it Proceedings of the VI ISPS Conference}, pp.3-8,1987;
{\it Indian Journal of Pure Applied Mathematics}, {\bf
22(11)}(1991), 887-903; {\it Journal of Multivariate Analysis}, {\bf
41(2)}(1992), 178-193; {\it Proceedings of the National  Academy of
Sciences}, {\bf LXV(II)}(1995), 121-142, {\bf LXV(III)}(1995),
227-251, {\bf LXVI(IV)}(1995), 367-393, {\bf LXVI(AI)}(1996), 1-22;
{\it Indian Journal of Pure and Applied Mathematics}, {\bf
24(9)}(1993), 513-531; {\it Advances in Applied Probability}, {\bf
31(2)}(1999), 343-354; {\it Rendiconti del Circolo Matematico di
Palermo, Series II, Suppl.}, {\bf 65}(2000), 219-232; {\it Linear
Algebra and Its Applications}, {\bf 183}(1993), 202-221; in {\it
Probability and Statistical Methods with Applications}, pp.293-316,
Chapman and Hall, 2001), Mathai and Saxena ({\it Journal de
Matematica e Estatistica}, {\bf 1}(1979), 91-106), Mathai and Rathie
({\it Statistica}, {\bf XL}(1980), 93-99; {\it Sankhya Series A}, {\bf
42}(1980), 78-87;), Mathai and Tracy ({\it Communications in
Statistics A}, {\bf 12(15)}(1983), 1727-1736; {\it Metron}, {\bf
44}(1986), 11-110), Mathai and Pederzoli ({\it Metron}, {\bf
LI(3-4)}(1993), 3-24; {\it Indian Journal of Pure Applied
Mathematics}, {\bf 27(3)}(1996), 7-32; {\it Linear Algebra and Its
Applications}, {\bf 253}(1997), 209-226, {\bf 269}(1998), 91-103). The important publication in this area is the book
on Jacobians of matrix transformation: A.M. Mathai, {\it
Jacobians of Matrix Transformations and Functions of Matrix
Argument}, World Scientific Publishing, New York, 1997. The work on
functions of matrix argument is continuing in the form of
applications in pathway models, fractional calculus and so on. These
will be mentioned later.\vskip.2cm In connection with matrix-variate
integrals it is a very often asked question that whether
matrix-variate integrals can be evaluated by treating them as
multiple integrals and by using standard techniques in calculus.
Mathai explored the possibility of explicitly evaluating
matrix-variate gamma and beta integrals as multiple integrals in
calculus. The basic matrix-variate integrals are the gamma integral
and beta integrals, where $X$ is a $p\times p$ real positive
definite matrix or hermitian positive definite matrix. For example,
when $X$ is real and $X>O$ (positive definite) the gamma integral is
$$\int_{X>O}|X|^{\alpha-\frac{p+1}{2}}{\rm e}^{-{\rm tr}(X)}{\rm d}X,\Re(\alpha)>\frac{p-1}{2}$$
and the beta integral is
$$\int_{O<X<I}|X|^{\alpha-\frac{p+1}{2}}|I-X|^{\beta-\frac{p+1}{2}}{\rm d}X,\Re(\alpha)>\frac{p-1}{2},\Re(\beta)>\frac{p-1}{2}.
$$The corresponding integrals are there in the complex-variate
case also. It is shown that this can be done explicitly for $p=2$
and a recurrence relation can be obtained so that step by step
they can be evaluated but for $p>2$ this method of treating as
multiple integrals is not a feasible proposition. See Mathai
({\it Journal of the Indian Mathematical Society}, {\bf 81(3-4)}(2014),
259-271; {\it Applied Mathematics and Computation}, {\bf 247}(2014), 312-318.)

\vskip.3cm\noindent{\bf 11.\hskip.3cm Multivariate Gamma and Beta Models}\\

\vskip.3cm Corresponding to a univariate model there is nothing
called a unique multivariate analogue. Explorations of some
convenient multivariate models corresponding to univariate gamma,
type-1 beta, type-2 beta, Dirichlet models etc were conducted in a
series of papers. See, for example, Mathai (In {\it Time Series
Methods in Hydrosciences,}, pp. 27-36, Elsevier, 1982), Mathai and
Moschopoulos ({\it Journal of Multivariate Analysis}, {\bf
39}(1991), 135-153; {\it Annals of the Institute of Statistical
Mathematics},{\bf 44(1)}(1992), 97-106; {\it Statistica}, {\bf
LVII(2)} (1992), 189-197, {\bf LIII(2)}(1993), 231-21). These were
some of the works on the multivariate analogues of gamma and beta
densities. Dirichlet models themselves are multivariate extensions
of type-1 and type-2 beta integrals or beta densities. When working
on order statistics from logistic populations, Mathai came across
the need for a generalized form of type-1 Dirichlet model, see
Mathai ({\it IEEE Trans. Reliability}, {\bf 52(2)}(2003), 20-206; in
{\it Statistical Methods and Practice: Recent Advances}, pp. 57-67,
Narosa Publishing, India, 2003; {\it Proceedings of the 7th
Conference of the Society for Special Functions and Their
Applications}, {\bf 7}(2006), pp. 131-142,). Various types of
generalizations of type-1 and type-2 Dirichlet densities were
considered, see for example, Jacob, Jose and Mathai ({\it Journal of
the Indian Academy of Mathematics}, {\bf 26(1)}(2004), 175-189);
Kurian, Kurian and Mathai ({\it Proceedings of the National Academy
of Sciences}, {\bf 74(A)II}(2004), 1-10), Jacob, George and Mathai
({\it Proceedings of the National  Academy of Sciences}, {\bf
15(3)}(2005, 1-9), Thomas and Mathai ({\it Advances in Applied
Statistics}, {\bf 8(1)}(2008), 37-56; {\it Sankhya Series A}, {\bf
71(1)}(2009), 49-63),
  Thomas, Thannippara and Mathai ({\it Journal of Probability and Statistical
  Science}, {\bf 6(2)}(2008), 187-200).\\

  \vskip.3cm\noindent{\bf 11.1.\hskip.3cm Power transformation and
  exponentiation}\\
  \vskip.3cm Another problem explored is to see the nature of models
  available by power transformations and exponentiation of
  standard probability models. Such a study is useful when
  looking for an appropriate model for a given data. These
  explorations are done in Mathai ({\it Journal of the Society
  for Probability and Statistics (ISPS)}, {\bf 13}(2012), 1-19).\\

  \vskip.3cm\noindent{\bf 11.2.\hskip.3cm Symmetric and asymmetric
  models}\\
  \vskip.3cm  A symmetric model,
  symmetric at $x=a$ where $a$ could be zero also, means that for
  $x<a$ the behavior of the function or the shape of the function
  is the same as its behavior for
  $x>a$. In many practical situations, symmetry may not be there.
  The behavior for $x<a$ may be different from that for $x>a$. Many
  authors have considered asymmetric models where asymmetry is
  introduced by giving different weighting factors for $x<a$ and for
  $x>a$ so that the total probability under the curve will be 1. But
  the shape of the curve itself may change for $x<a$ and for $x>a$.
  A method is proposed in the paper referred to in 11.1 above (Mathai 2012)
  where asymmetry is introduced through a scaling parameter so that
  the shape itself will be different for $x<a$ and $x>a$ cases but
  the total probability remaining as 1, which may have more practical relevance.

  \vskip.3cm\noindent{\bf 12.\hskip.3cm The Pathway Model}\\
  \vskip.3cm The basic idea was there in a paper of 1970's in the area
  of population studies where it was shown that by a limiting process
  one can go from one class of functions to another class of functions,
  the property is basically coming from the theory of hypergeometric
  functions from the aspect of getting rid off a numerator or a
  denominator parameter. This idea was revived and written as a
  paper on functions of matrix argument where the variable matrix
  is a rectangular one, see Mathai ({\it Linear Algebra and
  Its Applications}, {\bf 396}(2005), 317-328). Let $X$ be a real
  $m\times n$ matrix, $m\le n$ and of rank $m$ be a matrix variable.
  Let $A$ be $m\times m$ and $B$ be $n\times n$ constant nonsingular
  matrices. Consider the function

  $$f(X)=C~|AXBX'|^{\gamma}|I-(1-\alpha)AXBX'|^{\frac{\eta}{1-\alpha}},\eta>0\eqno(12.1)
  $$where $\alpha,\eta, C$ be scalar constants. This $C$ can act as
  a normalizing constant if we wish to create statistical density out
  of (12.1). Consider the case when $m=1,n=1$ and $x>0$. Then one
  can also take powers for $x$ and the model in (12.1) can be written as
  $$f_1(x)=c_1~x^{\gamma}[1-a(1-\alpha)x^{\delta}]^{\frac{\eta}{1-\alpha}}\eqno(12.2)
  $$where $a>0,\delta>0,\eta>0,x\ge 0$. In the matrix-variate case
  in (12.1) arbitrary powers for matrices is not feasible even
  though $AXBX'$ is positive definite because even for a positive
  definite matrix, $Y$, arbitrary power such as $Y^{\delta}$ may
  not be uniquely defined. Even when uniquely defined transformation
  such as $Z=Y^{\delta}$ will create problems when computing the Jacobians.
  The types of difficulties that can arise may be seen for the case
  $\delta=2$ described in the book, A.M. Mathai, {\it Jacobians of
  Matrix Transformations and Functions of Matrix Argument}, World
  Scientific Publishing, New York 1997. Hence for the matrix case we
  consider only when $\delta =1$. Consider case $-\infty<\alpha<1$.
  Then (12.2) remains as it is given in (12.2) which is a generalized
  type-1 beta function. But if $\alpha>1$ then writing $1-\alpha=-(\alpha-1)$
  the form in (12.2) changes to the following:

  $$f(x)=c_2~x^{\gamma}[1+a(\alpha-1)x^{\delta}]^{-\frac{\eta}{\alpha-1}}\eqno(12.3)
  $$for $a>0,\alpha>1,\eta>0,\delta>0,x\ge 0$. This model is a
  generalized type-2 beta model. When $\alpha\to 1$ in (12.2)
  and (12.3), $f_1(x)$ and $f_2(x)$ reduce to the the form
  $$f_3(x)=c_3~x^{\gamma}{\rm e}^{-a\eta x^{\delta}},a>0,\eta>o,x\ge 0.\eqno(12.4)
  $$This is a generalized gamma model. Thus three functional
  forms $f_1(x),f_2(x),f_3(x)$ are available for $\alpha<1.\alpha>1,\alpha\to 1$.
  This parameter $\alpha$ is called the pathway parameter,
  a pathway showing three different families of functions.
  
\vskip.2cm The practical utility of the model is that if (12.4)
  is the stable or ideal situation in a physical system then the
  unstable neighborhoods or functions leading to (12.4) are given
  in (12.2) and (12.3). In a model building situation, if the
  underlying data show a gamma-type behavior then a best-fitting
  model can be constructed for some values of the parameters or
  for some value of $\alpha$ the ideal model can be determined.
  Most of the statistical models in practical use in the areas
  of statistics, physics and engineering fields can be seen to
  be a member or products of members from $f_1,f_2,f_3$ above.
  Note that for $\alpha>1$ and $\alpha\to 1$ situations we can
  take $\delta>0$ or $\delta<0$ and both these situations can
  create statistical densities. Note that $f_1$ is a family of
  finite range models whereas $f_2$ and $f_3$ are families of
  infinite range models. Extended models are available by replacing
  $x$ by $|x|$ so that the whole real line will be covered. In
  this case the nonzero part of model (12.2) will be in the
  range $\pm [a(1-\alpha)]^{-\frac{1}{\delta}}$ and for others
  $-\infty<x<\infty$. Note that in (12.1) all individual variables
  $x_{ij}$'s are allowed to vary over the whole real line subject
  to the condition $I-(1-\alpha)AXBX'>O$ (positive definite). This
  model is also extended to complex rectangular matrix-variate case,
   see Mathai and Provost ({\it Linear Algebra and Its Applications},
   {\bf 410}(2005), 198-216).

  \vskip.2cm Note that (12.2) for $\gamma=0,\delta=1,a=1,\eta=1$ is
  Tsallis statistics in nonextensive statistical mechanics. The
  function, without the normalizing constant $c_1$ will then be
  $$g(x)=[1-(1-\alpha)x]^{\frac{1}{1-\alpha}}\eqno(12.5)
  $$
  which is Tsallis statistics. This can be generated by optimizing
  Tsallis entropy or Havrda-Charv\'at entropy with the denominator
  factor $1-\alpha$ instead of $2^{1-\alpha}-1$, subject the constraint
  that the first moment is fixed and this condition can be connected
  to the principle of the total energy being conserved. Note that (12.5) is also a power function model.
  $$\frac{{\rm d}}{{\rm d}x}g(x)=-[g(x)]^{\alpha}.
  $$Also (12.4) for $a=1,\delta=1,\eta=1$ is superstatistics in
  nonextensive statistical mechanics.
  \vskip.2cm Mathai's students have introduced a pathway fractional
  integral operator based on (12.2) and a pathway transform based on
  (12.2) and (12.3). (12.2),(12.3) can also be
  obtained by optimizing Mathai's entropy
  $$M_{\alpha}(f)=\frac{\int_{-\infty}^{\infty}[f(x)]^{2-\alpha}{\rm d}x-1}{\alpha-1},\alpha\ne 1,\alpha<2
  $$subject to two moment type constraints and also the pathway
  parameter $\alpha$ can be derived in terms of moments of $f_1(x)$
  or $f_2(x)$. Thus, in terms of entropies one can establish a
  entropic pathway, in terms of distributions as explained above one
  can create a distributional pathway, one can also look into the
  corresponding differential equations and create a differential
  pathway, covering the three sets of functions belonging to generalized
  and extended type-1 beta  family, type-2 beta family and gamma family.
  The theory of quadratic and bilinear forms in random variables is
  extended to cover pathway populations, instead of Gaussian population.
  Note that Gaussian population is a special case of the extended pathway
  population or pathway model, see Mathai$^{*}$ ({\it Linear Algebra and
  Its Applications}, {\bf 425}(2007), 162-170). Applications
  and advancement of theory of pathway model by Mathai and his associates
  may be seen from the following: Mathai and Haubold ({\it Physica A},
  {\bf 375}(2007), 110-122, {\bf 387}(2007), 2462-2470; {\it Physics
  Letters A}, {\bf 372}(2008), 2109-2113; {\it Integral Transforms and
  Special Functions}, {\bf 21(11)}(2011), 867-875; {\it Applied Mathematics
  and Computations}, {\bf 218}(2011), 799-804; {\it Mathematica Aeterna},
  {\bf 2(1)},(2012), 51-61; {\it Sun and Geosphere}, {\bf 8(2)}(2013), 63-70,
  {\it UN Proceedings (2013)};
  {\it Entropy}, {\bf 15}(2013), 4011-4025), Mathai and Provost ({\it IEEE
  Transactions on Reliability}, {\bf 55(2)}(2006), 237-244; {\it Journal
  of Probability and Statistical Science}, {\bf 9(1)}(2011), 1-20;
  {\it Physica A}, {\bf 392(4)}(2013), 545-551).\\

  \vskip.3cm\noindent{\bf 12.1.\hskip.3cm Input-output models}\\
  \vskip.3cm Many practical situations are input-output situations
  where what is observed is really the residual effect. Energy may
  be produced and consumed and what is observed is the net result or
  the residual effect. Water flows into a dam, which is the input
  variable,  and water is taken out
  of the dam, which is the output variable and the storage at any
  instant is the residual effect of the input minus the output. In
  any production-consumption, creation-destruction, growth-decay
  situation what is observed is $z=x-y$ where $x$ is the input
  variable and $y$ is the output variable and $z$ is the residual
  effect. Mathai explored a number of situations where $x$ and $y$
  are independently distributed real scalar random variables or
  matrix random variables. Observations as widely different as solar neutrinos and the amount
  of melatonin present in human body are all residual observations.
  Some works in this direction may be seen from Mathai ({\it Annals
  of the Institute of Statistical Mathematics}, {\bf 34}(1982),
  591-597; In {\it Time Series Methods in Hydrosciences}, pp. 27-36,
  Elsevier, Amsterdam, 1982; {\it Canadian Journal of Statistics},
  {\bf 21(3)}(1993), 277-283; {\it Journal of Statistical Research},
  {\bf 27(1-2)}(1993), 57-80; {\it Integral Transforms and Special Functions},
  {\bf 20(12)}(2009), 49-63), Haubold and Mathai ({\it Astrophysics
  and Space Science}, {\bf 228}(1995), 113-134; {\it Astrophysics and
  Space Science}, {\bf 273}(2000), 53-63), Saxena, Mathai and
  Haubold ({\it Astrophysics and Space Science}, {\bf 290}(2004),
  299-310) and a number of papers on fractional reaction-diffusion
  equations.\\

  \vskip.3cm\noindent{\bf 13.\hskip.3cm Work on Mittag-Leffler functions and Mittag-Leffler Densities}\\

  \vskip.3cm On Mittag-Leffler functions and their generalizations an overview
  paper is written, see Haubold, Mathai and Saxena ({\it Journal of Applied
  Mathematics}, {\bf ID 298628}(2011), 51 pages). Mittag-Leffler function comes
  in naturally when looking for solutions of fractional differential equations.
  This aspect will be considered later. Three standard forms of Mittag-Leffler
  functions in current use are the following:
  \begin{align*}
  E_{\alpha}(x)&=\sum_{k=0}^{\infty}\frac{x^k}{\Gamma(1+\alpha k)},\Re(\alpha)>0\\
  E_{\alpha,\beta}(x)&=\sum_{k=0}^{\infty}\frac{x^k}{\Gamma(\beta+\alpha k)},\Re(\alpha)>0,\Re(\beta)>0\\
  E_{\alpha,\beta}^{\gamma}(x)&=\sum_{k=0}^{\infty}\frac{(\gamma)_kx^k}{k!\Gamma(\beta+\alpha k)},\Re(\alpha)>0,\Re(\beta)>0.
  \end{align*}There is no condition on the parameter $\gamma$. If these are to be
  written in terms of H-functions then $\alpha$ and $\gamma$ have to be real and
  positive. A generalization can be made by introducing a general hypergeometric
  type function, which may be written as
  $$
  E_{\alpha,\beta,b_1,...,b_s}^{a_1,...,a_r}(x^{\delta})
  =\sum_{k=0}^{\infty}\frac{(a_1)_k...(a_r)_k(x^{\delta})^k}{k!\Gamma(\beta+\alpha k)(b_1)_k...(b_s)_k}
  $$where the notation $(a_j)_k$ and $(b_j)_k$ are Pochhammer symbols.
  Convergence conditions can be worked out for this general form.
  \vskip.2cm
  A problem of interest in this case is a general Mittag-Leffler density because
  such a density is needed in non-Gaussian stochastic processes and time series areas.
  Such a density was introduced based on $E_{\alpha,\beta}^{\gamma}(x^{\delta})$
  and it is shown that such a model is connected to fat-tailed models, L\'evy,
  Linnik models. Structural properties and asymptotic behavior are also
  studied and it is shown that such models are not attracted to Gaussian
  models, see Mathai ({\it Fractional Calculus \& Applied Analysis}, {\bf 13(1)}
  (2010), 113-132), Mathai and Haubold ({\it Integral Transforms and Special
  Functions}, {\bf 21(11)}(2011), 867-875).\\

  \vskip.3cm\noindent{\bf 14.\hskip.3cm Work on Kr\"atzel Function and Kr\"atzel Densities}\\
  \vskip.3cm Another area explored is Kr\"atzel function, Kr\"atzel transform and Kr\"atzel densities. Since Kr\"atzel transform is important in applied analysis area, a general density is introduced based on Kr\"atzel integral. The basic Kr\"atzel integral is of the form
$$g_1(x)=\int_0^{\infty}x^{\gamma}{\rm e}^{-ax-\frac{y}{x}}{\rm d}x,a>0,y>0\eqno(14.1)
$$which can be generalized to the form
$$g_2(x)=\int_0^{\infty}x^{\gamma}{\rm e}^{-ax^{\alpha}-\frac{y}{x^{\beta}}}{\rm d}x\eqno(14.2)
$$for $a>0,y>0,\alpha>0,\beta>0$ or $\beta<0$. The integrand in (14.1), normalized, is the inverse Gaussian density. The integral itself can be interpreted as Mellin convolution of a product, the marginal density in a bivariate case etc. The integral in (14.2) is connected the general reaction-rate probability integral in reaction-rate theory ($\beta=\frac{1}{2},\alpha=1$ is the basic integral in reaction-rate theory) , unconditional densities in Bayesian analysis, marginal densities in a bivariate set up, and so on. Different problems in a large number of areas can be connected to (14.2). Note that $x^{\gamma}{\rm e}^{-ax^{\alpha}}$, normalized can act as a marginal density of a real scalar random variable $x>0$ and ${\rm e}^{-\frac{y}{x^{\beta}}}$, normalized, can act as the conditional density of $y$, given $x$. In this case (14.2) has the structure of unconditional density of $y$ in a Bayesian analysis situation. One can also look at (14.2), normalized, as the joint density of two real scalar positive random variables and in this case (14.2) integral represents the marginal density of $y$. For $\beta>0$, (14.2) can act as the Mellin convolution of a product and for $\beta<0$ it can represent the Mellin convolution of a ratio. In this case, one can connect it to the Laplace transform of a generalized gamma density for $\alpha=1$. Many types of such properties are studied in Mathai ({\it International Journal of Mathematical Analysis}, {\bf 6(51)}(2012), 2501-2510; In {\it Frontiers of Statistics and its Applications}, Bonfring Publications, Germany, 2013;
{\it Proceedings of the 10th and 11th Annual Conference of SSFA}, {\bf 10-11}(2011-2012), pp. 11-20). Mathai has also considered the matrix-variate version of (14.1).

  \vskip.3cm\noindent{\bf 15.\hskip.3cm Work on Fractional Calculus}\\

  \vskip.3cm Mathai may be credited with making a connection of fractional integrals
  to statistical distribution theory, extending fractional calculus to matrix-variate cases,
  to complex matrix-variate cases, to many scalar variable (multiple) cases,
  to many matrix variable cases. Recently Mathai has given a geometrical interpretation
  of fractional integrals in a simplex as fractions of certain total integral in
  $n$-dimensional cube. Mathai has also given a new definition to the area of fractional
  integrals, and thereby fractional derivatives, as Mellin convolutions of products and
  ratios in the real scalar case and as M-convolutions of products and ratios in the
  matrix-variate case, where one function is of type-1 beta form, see Mathai ({\it Integral Transforms and Special Functions},
  {\bf 20(12)}(2009), 871-882; {\it Linear Algebra and Its Applications}, {\bf 439}(2013),
  2901-2913, {\bf 446}(2014), 196-215), Mathai and
  Haubold ({\it Fractional Calculus \& Applied Analysis}, {\bf 14(1)}(2011), 138-155;
  {\it Cornell University arXiv}, {I-IV}(2012) 4 papers; {\it Fractional Calculus
  \& Applied Analysis}, {\bf 16(2)}(2013), 469-478). Papers are
  published solving various types of fractional reaction, diffusion, reaction-diffusion
  differential equations, see Haubold, Mathai, Saxena ({\it Bulletin of the Astronomical
  Society of India}, {\bf 35}(2007), 681-689; {\it Journal of Computational and Applied
  Mathematics}, {\bf 235}(2011), 1311-1316; {\it Journal of Mathematical Physics},
  {\bf 51}(2010), 103506-8), Saxena, Mathai and Haubold ({\it Astrophysics and
  Space Science}, {\bf 305}(2006), 289-296, 297-303, 305-313; {\it Astrophysics
  and Space Science Proceedings}, (2010), pp.35-40, 55-62; {\it Axiom}, {\bf 3(3)}
  (2014), 320-334; {\it Journal of Mathematical Physics}, {\bf 55}(2014), 083519
   doi:10.1063/1.4891922.

\end{document}